\documentclass[graybox]{svmult}
\usepackage{amsmath,amssymb,amsthm,mathtools}

\usepackage{type1cm}        
\usepackage{makeidx}         
\usepackage{graphicx}        
\usepackage{multicol}        
\usepackage[bottom]{footmisc}

\usepackage{newtxtext}       %
\usepackage[varvw]{newtxmath}       

\makeindex             

\begin{document}

\title*{Fractional Sobolev spaces via interpolation, and applications to mixed local-nonlocal operators}
\titlerunning{Fractional Sobolev spaces, interpolation, and mixed local-nonlocal operators}
\author{Alberto Maione\orcidID{0000-0003-1629-6670}}
\authorrunning{A. Maione}
\institute{Alberto Maione \at Centre de Recerca Matemàtica, Campus de Bellaterra, Edifici C, 08193 Bellaterra, Spain, \newline\email{amaione@crm.cat}}
%
%
\maketitle

\abstract*{In this note, we present a well-known connection between the Sobolev-Slobodeckij spaces, also known as \emph{Fractional Sobolev spaces}, and \emph{interpolation theory}.
We show how Sobolev spaces can be equivalently characterized as real and complex interpolation spaces between Lebesgue spaces and integer-order Sobolev spaces.
We also state a spectral theorem for the so-called \emph{mixed local-nonlocal} operators, and show how interpolation theory leads to its proof.
This note is intended for early-career researchers, and aims to provide a concise and accessible introduction to the subject.}

\abstract{In this note, we present a well-known connection between the Sobolev-Slobodeckij spaces, also known as \emph{Fractional Sobolev spaces}, and \emph{interpolation theory}.
We show how Sobolev spaces can be equivalently characterized as real and complex interpolation spaces between Lebesgue spaces and integer-order Sobolev spaces.
We also state a spectral theorem for the so-called \emph{mixed local-nonlocal} operators, and show how interpolation theory leads to its proof.
This note is intended for early-career researchers, and aims to provide a concise and accessible introduction to the subject.}


\section{Interpolation theory}
We recall some classical definitions and results from interpolation theory in the setting of \emph{Banach spaces}.
We refer the interested reader to \cite{A,BL,G,Leo,L,S,T} for a more detailed exposition of the topic.

\begin{definition}
Let $(X,\|\cdot\|_X)$ and $(Y,\|\cdot\|_Y)$ be Banach spaces.
We say that the couple $(X,Y)$ is \emph{compatible} if there exists a Hausdorff topological vector space $V$ such that $X\subseteq V$ and $Y\subseteq V$.
\end{definition}

Given a compatible couple of Banach spaces $(X,Y)$, we can define the \emph{sum} and \emph{intersection} spaces as follows
\begin{align*}
    X+Y&:=\{f\in V:f=g+h\text{ with }g\in X\text{ and }h\in Y\}\subseteq V,\\
    X\cap Y&:=\{f\in V: f\in X\text{ and }f\in Y\}\subseteq V.
\end{align*}
By definition, it holds that
\[
X\cap Y\subseteq X\cup Y\subseteq X+Y.
\]
We work with the sum space $X+Y$ instead of the union space $X\cup Y$, since the latter is not a vector space, even though both $X$ and $Y$ are Banach spaces.
We associate $X+Y$ and $X\cap Y$ with the graph norms, respectively given by
\begin{align}
    \|f\|_{X+Y}&:=\inf_{f=g+h, g\in X, h\in Y}\{\|g\|_X+\|h\|_Y\},\label{norms}\\
    \|f\|_{X\cap Y}&:=\max(\|f\|_X,\|f\|_Y).\notag
\end{align}

\begin{theorem}
Assume that $(X,\|\cdot\|_X)$ and $(Y,\|\cdot\|_Y)$ are normed spaces.
Then, the spaces $(X+Y,\|\cdot\|_{X+Y})$ and $(X\cap Y,\|\cdot\|_{X\cap Y})$ are normed spaces.
Moreover, if $X,Y$ are complete, then $X+Y$ and $ X\cap Y$ are complete as well.
\end{theorem}

\begin{definition}
Let $(X,Y)$ be a compatible couple of Banach spaces and let $W$ be a Banach space.
We say that $W$ is an \emph{intermediate space} between $X$ and $Y$ if
\[
X\cap Y\subseteq W\subseteq X+Y.
\]
All inclusions are continuous.
\end{definition}

A fundamental example comes from classical functional analysis.
Let $\Omega$ be an open set of finite measure and let $p\leq r\leq q$. Consider the following Banach spaces
\begin{align*}
    X=L^p(\Omega),\quad Y=L^q(\Omega),\quad W=L^r(\Omega).
\end{align*}
By the nesting property of Lebesgue spaces, $L^q(\Omega)\hookrightarrow L^r(\Omega)\hookrightarrow L^p(\Omega)$ and, up to equivalent norms,
\[
X+Y=L^p(\Omega),\quad X\cap Y=L^q(\Omega).
\]
Then, $W=L^r(\Omega)$ is an intermediate space between $X=L^p(\Omega)$ and $Y=L^q(\Omega)$, since
\[
X\cap Y=L^q(\Omega)\subseteq W=L^r(\Omega)\subseteq L^p(\Omega)=X+Y.
\]

The space $L^r(\Omega)$ of the previous example is not only an intermediate space between $L^p(\Omega)$ and $L^q(\Omega)$, but actually an \emph{interpolation space}, according to the following definition.
We use the following notation.
\begin{itemize}
    \item $\mathcal{L}(X)$ is the set of linear and bounded transformations $T:X\to X$.
    \item $\mathcal{L}(X)\cap\mathcal{L}(Y)$ is the set of $T:X+Y\to X+Y$ such that $T|_X\in\mathcal{L}(X)$ and $T|_Y\in\mathcal{L}(Y)$.
    \item $\mathcal{L}(X_0,X_1)$ is the set of linear and bounded transformations $T:X_0\to X_1$.
    \item $\mathcal{L}(X_0,X_1)\cap\mathcal{L}(Y_0,Y_1)$ is the set of $T:X_0+Y_0\to X_1+Y_1$ such that $T|_{X_0}\in\mathcal{L}(X_0,X_1)$ and $T|_{Y_0}\in\mathcal{L}(Y_0,Y_1)$.
\end{itemize}

\begin{definition}
Let $(X,Y)$ be a compatible couple of Banach spaces and let $W$ be an intermediate space between $X$ and $Y$.
We say that $W$ is an \emph{interpolation space} between $X$ and $Y$ if for every $T\in\mathcal{L}(X)\cap\mathcal{L}(Y)$ we have
\[
T|_W\in\mathcal{L}(W).
\]
\end{definition}

The reason why, in the previous example, the space $W=L^r(\Omega)$ is an interpolation space between $X=L^p(\Omega)$ and $Y=L^q(\Omega)$ directly follows from the following lines.

\begin{theorem}[M. Riesz, 1926; Thorin, 1938, 1948]\label{Thm:Riesz-Thorin}
Let $1\leq p_0\leq p_1\leq+\infty$, and let $1\leq q_0\leq q_1\leq+\infty$, with $p_0\leq q_0$ and $p_1\leq q_1$.\footnote{As proved by Thorin in 1938, the assumption $p_0\leq q_0$ and $p_1\leq q_1$ may be dropped.}
Let $T\in\mathcal{L}(L^{p_0},L^{q_0})\cap\mathcal{L}(L^{p_1},L^{q_1})$.
For $s\in(0,1)$, let $p_s\in[p_0,p_1]$, with $\frac{1}{p_s}=(1-s)\frac{1}{p_1}+s\frac{1}{p_0}$, and let $q_s\in[q_0,q_1]$, with  $\frac{1}{q_s}=(1-s)\frac{1}{q_1}+s\frac{1}{q_0}$.
Then, $T\in\mathcal{L}(L^{p_s},L^{q_s})$ and
\begin{equation}\label{Interpolation_inequality_Riesz_Thorin}
\|T\|_{\mathcal{L}(L^{p_s},L^{q_s})}:=\sup_{h\in L^{p_s}, h\neq 0}\frac{\|Th\|_{L^{q_s}}}{\|h\|_{L^{p_s}}}\leq \|T\|_{\mathcal{L}(L^{p_1},L^{q_1})}^{1-s}\|T\|_{\mathcal{L}(L^{p_0},L^{q_0})}^s,
\end{equation}
where
\[
\|T\|_{\mathcal{L}(L^{p_0},L^{q_0})}:=\sup_{f\in L^{p_0}, f\neq 0}\frac{\|Tf\|_{L^{q_0}}}{\|f\|_{L^{p_0}}},\quad
\|T\|_{\mathcal{L}(L^{p_1},L^{q_1})}:=\sup_{f\in L^{p_1}, f\neq 0}\frac{\|Tf\|_{L^{q_1}}}{\|f\|_{L^{p_1}}}.
\]
\end{theorem}

\begin{definition}
Let $(X_0,Y_0)$ and $(X_1,Y_1)$ be compatible couples of Banach spaces, and let $W_0$ and $W_1$ be Banach spaces.
We say that $(W_0,W_1)$ is an \emph{interpolation couple} with respect to $(X_0,Y_0)$ and $(X_1,Y_1)$ if:
\begin{itemize}
    \item[$(a)$]\, $W_0$ is an intermediate space between $X_0$ and $Y_0$.
    \item[$(b)$]\, $W_1$ is an intermediate space between $X_1$ and $Y_1$.
    \item[$(c)$]\, For every $T\in\mathcal{L}(X_0,X_1)\cap\mathcal{L}(Y_0,Y_1)$ it holds that $T\in\mathcal{L}(W_0,W_1)$.
\end{itemize}
\end{definition}

\begin{remark}
Let $(W_0,W_1)$ be an interpolation couple with respect to $(X_0,Y_0)$ and $(X_1,Y_1)$.
Then, in general, $W_0$ is \emph{not} an interpolation space between $X_0$ and $Y_0$ (equivalently, $W_1$ is not an interpolation space between $X_1$ and $Y_1$).
However, in practical situations, it is convenient to work with couples of spaces, since by letting $X:=X_0=X_1$ and $Y:=Y_0=Y_1$, then $W_0\equiv W_1=:W$, and
\begin{center}
    $W$ becomes an interpolation space between $X$ and $Y$.
\end{center}
\end{remark}

By means of the previous definition, we can reinterpret Theorem~\ref{Thm:Riesz-Thorin} as follows.
For $1\leq p_0\leq p_1\leq+\infty$ and $1\leq q_0\leq q_1\leq+\infty$, let $(L^{p_0},L^{p_1})$ and $(L^{q_0},L^{q_1})$ be compatible couples.
For $s\in(0,1)$, let $r_i\in[p_i,q_i]$, with
\[
\frac{1}{r_i}=(1-s)\frac{1}{q_i}+s\frac{1}{p_i}\quad (i=0,1).
\]
Then, $(L^{r_0},L^{r_1})$ is an interpolation couple between the couples $(L^{p_0},L^{p_1})$ and $(L^{q_0},L^{q_1})$.
Additionally, if $p_0=q_0=:p$, $p_1=q_1=:q$, and $p_s=q_s=:r$, we get that $L^r$ is an interpolation space between $L^p$ and $L^q$, and we recover the following interpolation inequality \eqref{Interpolation_inequality}, as a consequence of \eqref{Interpolation_inequality_Riesz_Thorin}.

\begin{theorem}[Interpolation inequality]\label{thm:interpolation_ineq}
Let $\Omega$ be a bounded open subset of $\mathbb{R}^n$, let $1\leq p\leq q\leq +\infty$, and let $f\in L^p(\Omega)\cap L^q(\Omega)$.
Then, $f\in L^r(\Omega)$ for every $r\in[p,q]$.
Moreover, the following estimate holds
\begin{equation}\label{Interpolation_inequality}
    \|f\|_{L^r(\Omega)}\leq\|f\|_{L^q(\Omega)}^{1-s}\|f\|_{L^p(\Omega)}^s,
\end{equation}
for some $s\in[0,1]$ such that $\frac{1}{r}=(1-s)\frac{1}{q}+s\frac{1}{p}$.
\end{theorem}


Other trivial examples of interpolation couples are as follows.
Let $(X_0,Y_0)$ and $(X_1,Y_1)$ be compatible couples of Banach spaces.
Denoted by
\begin{align*}
    W_0=X_0+Y_0&\text{ and }W_1=X_1+Y_1,\text{ and}\\
    W_0=X_0\cap Y_0&\text{ and }W_1=X_1\cap Y_1,
\end{align*}
then, the couples $(W_0,W_1)$ are interpolation couples with respect to $(X_0,Y_0)$ and $(X_1,Y_1)$.

Let $(X_0,Y_0)$ and $(X_1,Y_1)$ be two compatible couples of Banach spaces, and consider $\mathcal{L}(X_0+Y_0,X_1+Y_1)$ and $\mathcal{L}(X_0\cap Y_0,X_1\cap Y_1)$, the sets of linear and bounded transformations
\begin{align*}
    T:X_0+Y_0\to X_1+Y_1,\quad T:X_0\cap Y_0\to X_1\cap Y_1.
\end{align*}
The associated norms $\|T\|_{\mathcal{L}(X_0+Y_0,X_1+Y_1)}$ and $\|T\|_{\mathcal{L}(X_0\cap Y_0,X_1\cap Y_1)}$ both satisfy
\begin{align*}
    \|T\|_{\mathcal{L}(X_0+Y_0,X_1+Y_1)}&:=\sup_{f\in X_0+Y_0, f\neq 0}\frac{\|Tf\|_{X_1+Y_1}}{\|f\|_{X_0+Y_0}}\leq\max\Big(\|T\|_{\mathcal{L}(X_0,X_1)},\|T\|_{\mathcal{L}(Y_0,Y_1)}\Big)\\
    \|T\|_{\mathcal{L}(X_0\cap Y_0,X_1\cap Y_1)}&:=\sup_{f\in X_0\cap Y_0, f\neq 0}\frac{\|Tf\|_{X_1\cap Y_1}}{\|f\|_{X_0\cap Y_0}}\leq\max\Big(\|T\|_{\mathcal{L}(X_0,X_1)},\|T\|_{\mathcal{L}(Y_0,Y_1)}\Big).
\end{align*}

\begin{definition}
Let $(X_0,Y_0)$ and $(X_1,Y_1)$ be compatible couples of Banach spaces, and let $(W_0,W_1)$ be an interpolation couple with respect to $(X_0,Y_0)$ and $(X_1,Y_1)$.
We say that the interpolation is \emph{exact} if for every $T\in\mathcal{L}(X_0,X_1)\cap\mathcal{L}(Y_0,Y_1)$, which gives $T\in\mathcal{L}(W_0,W_1)$, it holds that
\[
\|T\|_{\mathcal{L}(W_0,W_1)}\leq\max\Big(\|T\|_{\mathcal{L}(X_0,X_1)},\|T\|_{\mathcal{L}(Y_0,Y_1)}\Big).
\]
\end{definition}

This requirement is quite restrictive. A weaker version is the following.

\begin{definition}
Let $(X_0,Y_0)$ and $(X_1,Y_1)$ be compatible couples of Banach spaces, and let $(W_0,W_1)$ be an interpolation couple with respect to $(X_0,Y_0)$ and $(X_1,Y_1)$.
We say that the interpolation is \emph{uniform}, if for every $T\in\mathcal{L}(X_0,X_1)\cap\mathcal{L}(Y_0,Y_1)$, there exists a positive constant $C$, only depending on $X_0,X_1,Y_0,Y_1$, such that
\[
\|T\|_{\mathcal{L}(W_0,W_1)}\leq C\max\Big(\|T\|_{\mathcal{L}(X_0,X_1)},\|T\|_{\mathcal{L}(Y_0,Y_1)}\Big).
\]
\end{definition}

In the setting of Banach spaces, it is shown in \cite[Theorem 2.4.2]{BL} that
\[
(W_0,W_1)\text{ is an interpolation couple if and only if the interpolation is uniform}.
\]
To the best of our knowledge, even in general normed space, there are no known counterexamples.
We note that the interpolation in Theorem~\ref{Thm:Riesz-Thorin} is exact of power $s$, according to the following definition.

\begin{definition}\label{def:exact_power_s}
Let $(X_0,Y_0)$ and $(X_1,Y_1)$ be compatible couples of Banach spaces, and let $(W_0,W_1)$ be an interpolation couple with respect to $(X_0,Y_0)$ and $(X_1,Y_1)$.

We say that the interpolation is of \emph{power} $s$, if for every $T\in\mathcal{L}(X_0,X_1)\cap\mathcal{L}(Y_0,Y_1)$, there exist $s\in(0,1)$ and a positive constant $C$, only depending on $X_0,X_1,Y_0,Y_1$, such that
\[
\|T\|_{\mathcal{L}(W_0,W_1)}\leq C\|T\|_{\mathcal{L}(X_0,X_1)}^{1-s}\|T\|_{\mathcal{L}(Y_0,Y_1)}^s.
\]
If $C=1$, we say that $(W_0,W_1)$ is also exact (of power $s$).
\end{definition}

We conclude this section with the following result, whose proof can be found e.g. in \cite[Theorem 2.5.1]{BL}.

\begin{theorem}[Aronszajn-Gagliardo, 1965]
Let $(X_0,Y_0)$ be a compatible couple of Banach spaces, and let $W_0$ be an interpolation space with respect to $(X_0,Y_0)$.
Then, there exists a compatible couple of Banach spaces $(X_1,Y_1)$, and a Banach space $W_1$ such that $W_1$ is an interpolation space with respect to $(X_1,Y_1)$, and $(W_0,W_1)$ is an exact interpolation couple with respect to $(X_0,Y_0)$ and $(X_1,Y_1)$.
\end{theorem}


\section{Interpolation methods}

\subsection{Real Interpolation. The $K$-method}
There are several methods to construct interpolation spaces.
In this section, we briefly recall the construction of \emph{real interpolation spaces} between Banach spaces through the so-called $K$-method by Peetre.
The construction proceeds in three steps.
\begin{itemize}
    \item [$1$] Given $(X,Y)$ a compatible couple of Banach spaces, we define an equivalent norm, with respect to \eqref{norms}, for the Banach space $X+Y$.
    \item [$2$] We define the real interpolation space $(X,Y)_{s,p}$, for $p\in[1,+\infty)$, and $s\in(0,1)$.
    \item [$3$] We show that $(X,Y)_{s,p}$ in an interpolation space between $X$ and $Y$.
\end{itemize}

\noindent [\textbf{Step 1:}] 
Let $(X,Y)$ be a compatible couple of Banach spaces and fix $x>0$.
We define the Peetre $K$-functional as
\begin{align*}
K: X+Y &\to \mathbb{R}^+_0 \\
f &\mapsto K(x,f,X,Y) := \inf_{f=g+h, g\in X, h\in Y}\Big\{\|g\|_X+x\|h\|_Y\Big\}.
\end{align*}
The $K$-functional defines a norm equivalent to $\|\cdot\|_{X+Y}$ (indeed, they coincide when $x=1$).\footnote{When clear from the context, we abbreviate $K(x,f,X,Y)$ with $K(x,f)$.}
We note that for every $x\neq 1$, the function is \emph{not symmetric} with respect to $X$ and $Y$, but instead satisfies the relation
\begin{equation}\label{non_symmetry}
K(x,f,Y,X)=xK(x^{-1},f,X,Y)\neq K(x,f,X,Y).
\end{equation}
This lack of symmetry on the $K$-functional is compensated by the symmetric property 
\[
(X,Y)_{s,p}=(Y,X)_{1-s,p}\quad\text{for every $s\in(0,1)$}.
\]

\noindent [\textbf{Step 2:}] 
We first recall the definition of weighted Lebesgue spaces.
\begin{definition}
Let $\Omega$ be an open subset of $\mathbb{R}^n$ and $p\in[1,+\infty)$, and let $\omega : \Omega \to (0,+\infty)$ be a measurable weight, positive a.e. 
We define the weighted Lebesgue space
\[
L^p_\omega(\Omega):=\left\{f\in L^p(\Omega):\int_\Omega|f(x)|^p\omega(x)\,\mathrm{d}x<\infty\right\},
\]
endowed with the norm
\[
\|f\|_{L^p_\omega(\Omega)}:=\left( \int_\Omega |f(x)|^p \, \omega(x)\, dx \right)^{1/p}.
\]
\end{definition}

\begin{definition}
Let $p\in[1,+\infty)$, $s\in(0,1)$ and let $x>0$.
The real interpolation space $(X,Y)_{s,p}$ is the subspace of $X+Y$ defined as
\[
(X,Y)_{s,p}:=\left\{f\in X+Y:\frac{K(x,f,X,Y)}{x^s}\in L^p_\frac{1}{x}(\mathbb{R}^+)\right\}.
\]
We endow $(X,Y)_{s,p}$ with the norm
\begin{align*}
    \|f\|_{(X,Y)_{s,p}}&=\left\|\frac{K(x,f,X,Y)}{x^s}\right\|_{L^p_\frac{1}{x}(\mathbb{R}^+)}\\
    &= \left( \int_0^{+\infty} \left| \frac{K(x,f,X,Y)}{x^s} \right|^p \frac{1}{x} \, dx \right)^{1/p} \\
    &= \left( \int_0^{+\infty} \frac{K(x,f,X,Y)^p}{x^{sp+1}} \, dx \right)^{1/p}.
\end{align*}
\end{definition}

\noindent [\textbf{Step 3:}] Properties of the real interpolation spaces.

For every $p\in[1,+\infty)$ and $s\in(0,1)$, $\big((X,Y)_{s,p},\|\cdot\|_{(X,Y)_{s,p}}\big)$ is a Banach space. 
Moreover, the following symmetry property follows by \eqref{non_symmetry}
\[
(X,Y)_{s,p}=(Y,X)_{1-s,p}.
\]

\begin{theorem}\label{thm:inclusions_interpolation}
Let $p_1,p_2\in[1,+\infty)$ and let $s_1,s_2\in(0,1)$. 
\begin{itemize}
    \item If $1\leq p_1\leq p_2\leq+\infty$, then for every $s\in(0,1)$
    \[
    X\cap Y\subseteq (X,Y)_{s,p_1}\subseteq (X,Y)_{s,p_2}\subseteq X+Y.
    \]
    \item If $Y\subseteq X$ and $s_1<s_2$ satisfy $s_2-\frac{1}{p_1}>s_1-\frac{1}{p_2}$, then
    \[
    (X,Y)_{s_2,p_1}\subseteq (X,Y)_{s_1,p_2}.
    \]
    In particular, for every $p\in[1,+\infty)$
    \begin{equation}\label{ext:Brezis_Miro}
         (X,Y)_{s_2,p}\subseteq (X,Y)_{s_1,p}.
    \end{equation}
\end{itemize}
All the above inclusions are continuous embeddings.
\end{theorem}

Let $p_1=p_2=p$. By Theorem~\ref{thm:inclusions_interpolation}, $(X,Y)_{s,p}$ is an intermediate space between $X,Y$.

\begin{theorem}\label{thm:real_interpolation}
Let $p\in[1,+\infty)$ and $s\in(0,1)$.
Let $(X_0,Y_0)$ and $(X_1,Y_1)$ be compatible couples of Banach spaces.
Let $T\in \mathcal{L}(X_0,X_1)\cap \mathcal{L}(Y_0,Y_1)$. Then $T$ extends to a bounded linear operator
\[
T\in\mathcal{L}((X_0,Y_0)_{s,p},(X_1,Y_1)_{s,p}).
\]
Moreover, for every $T\in\mathcal{L}(X_0,X_1)\cap\mathcal{L}(Y_0,Y_1)$, it holds that
\[
\|T\|_{\mathcal{L}((X_0,Y_0)_{s,p},(X_1,Y_1)_{s,p})}\leq\|T\|_{\mathcal{L}(X_0,X_1)}^{1-s}\|T\|_{\mathcal{L}(Y_0,Y_1)}^s.
\]
\end{theorem}

Let $X_0=X_1=:X$ and $Y_0=Y_1=:Y$.
By Theorem~\ref{thm:real_interpolation}, $(X,Y)_{s,p}$ is an interpolation space between $X$ and $Y$. The interpolation is exact of power $s$; see Definition~\ref{def:exact_power_s}.

\begin{corollary}
Let $(X,Y)$ be a compatible couple of Banach spaces. 
Then, for every $p\in[1,+\infty)$ and $s\in(0,1)$ there exists a constant $c$, depending only on $p$ and $s$, such that
\[
\|f\|_{(X,Y)_{s,p}}\leq c\|f\|_Y^{1-s}\|f\|_X^s\quad\text{for every }f\in X\cap Y.
\]
\end{corollary}

\subsection{Complex Interpolation} 
In this section, we briefly recall the construction of \emph{complex interpolation spaces} between Banach spaces.
\medskip

\noindent [\textbf{Step 1:}] Construction of the complex interpolation spaces.

Let $(X,Y)$ be a compatible couple of Banach spaces.
Consider the vector space $\mathcal{F}(X,Y)$ of functions $f:\mathbb{C}\to X+Y$ satisfying:
\begin{itemize}
    \item [$(a)$]\, $f$ is holomorphic (analytic) on the set $S=\{z\in\mathbb{C}: 0<\textrm{Re}(z)<1\}$,
    \item [$(b)$]\, $f$ is continuous and bounded on $\overline{S}=\{z\in\mathbb{C}: 0\leq \textrm{Re}(z)\leq 1\}$,
    \item [$(c)$]\, the mapping $t\mapsto f(it)$ is continuous from $\mathbb{R}$ to $X$,
    \item[$(d)$]\, the mapping $t\mapsto f(1+it)$ is continuous from $\mathbb{R}$ to $Y$. 
\end{itemize}
We endow $\mathcal{F}(X,Y)$ with the norm
\begin{align*}
    \|f\|_{\mathcal{F}(X,Y)}:=\max\Big\{\sup_{t\in\mathbb{R}}\|f(it)\|_X,\sup_{t\in\mathbb{R}}\|f(1+it)\|_Y\Big\}.
\end{align*}
Then, $(\mathcal{F}(X,Y),\|\cdot\|_{\mathcal{F}(X,Y)})$ is a Banach space.

\begin{definition}
Let $s\in(0,1)$. The complex interpolation space $(X,Y)_{[s]}$ is the subspace of $X+Y$ defined as
\[
(X,Y)_{[s]} := \left\{ g : \exists f \in \mathcal{F}(X,Y) \text{ such that } g(t) = f(s + it) \text{ for } t \in \mathbb{R} \right\}.
\]
We endow $(X,Y)_{[s]}$ with the norm
\[
\|g\|_{(X,Y)_{[s]}}:=\inf\big\{\|f\|_{\mathcal{F}(X,Y)}: f\in\mathcal{F}(X,Y), g=f(s+i\cdot)\big\}.
\]
\end{definition}

\noindent [\textbf{Step 2:}] Properties of the complex interpolation spaces.

For every $s\in(0,1)$, $\big((X,Y)_{[s]},\|\cdot\|_{(X,Y)_{[s]}}\big)$ is a Banach space. 
Moreover, the following symmetry property follows by definition
\[
(X,Y)_{[s]}=(Y,X)_{[1-s]}.
\]

\begin{theorem}\label{thm:inclusions_interpolation_complex}
Let $p\in[1,\infty)$ and let $s,s_1,s_2\in(0,1)$. 
\begin{itemize}
    \item If $Y\subseteq X$ and $s_1<s_2$, then
    \[
    (X,Y)_{[s_2]}\subseteq (X,Y)_{[s_1]}.
    \]
    \item For every $s\in(0,1)$, we have
    \[
    (X,Y)_{s,1}\subseteq (X,Y)_{[s]}\subseteq (X,Y)_{s,\infty}.
    \]
    \item In general, $(X,Y)_{s,p}\neq (X,Y)_{[s]}$ unless $p=2$ (or in special cases such as Hilbert couples).
\end{itemize}
All the above inclusions are continuous embeddings.
\end{theorem}

For every $s \in (0,1)$, the space $(X,Y)_{[s]}$ is an intermediate space between $X$ and $Y$.

\begin{theorem}\label{thm:complex_interpolation}
Let $p\in[1,+\infty)$ and $s\in(0,1)$.
Let $(X_0,Y_0)$ and $(X_1,Y_1)$ be compatible couples of Banach spaces.
Let $T\in \mathcal{L}(X_0,X_1)\cap \mathcal{L}(Y_0,Y_1)$. Then $T$ extends to a bounded linear operator
\[
T\in\mathcal{L}((X_0,Y_0)_{[s]},(X_1,Y_1)_{[s]}).
\]
Moreover, the operator norm satisfies the interpolation estimate
\[
\|T\|_{\mathcal{L}((X_0,Y_0)_{[s]},(X_1,Y_1)_{[s]})}\leq\|T\|_{\mathcal{L}(X_0,X_1)}^{1-s}\|T\|_{\mathcal{L}(Y_0,Y_1)}^s.
\]
\end{theorem}

Let $X_0=X_1=:X$ and $Y_0=Y_1=:Y$.
By Theorem~\ref{thm:complex_interpolation}, $(X,Y)_{[s]}$ is an interpolation space between $X$ and $Y$. The interpolation is exact of power $s$; see Definition~\ref{def:exact_power_s}.

\begin{corollary}
Let $(X,Y)$ be a compatible couple of Banach spaces. 
For every $s\in(0,1)$ there exists a constant $c$, depending only on $s$, such that
\[
\|f\|_{(X,Y)_{[s]}}\leq c\|f\|_Y^{1-s}\|f\|_X^s\quad\text{for every }f\in X\cap Y.
\]
\end{corollary}


\section{Sobolev spaces via interpolation}

\begin{definition}
Let $\Omega$ be a bounded open subset of $\mathbb{R}^n$, let $p\in[1,+\infty)$ and let $s\in\mathbb{R}^+$.
\begin{itemize}
    \item We define the Sobolev space $W^{1,p}(\Omega)$ as the space
    \[
    W^{1,p}(\Omega):=\Big\{f\in L^p(\Omega): \nabla f\in L^p(\Omega;\mathbb{R}^n)\Big\}
    \]
    and we endow $W^{1,p}(\Omega)$ with the graph norm
    \begin{align*}
        \|f\|_{W^{1,p}(\Omega)}:=    \left(\|f\|^p_{L^p(\Omega)}+\|\nabla f\|^p_{L^p(\Omega;\mathbb{R}^n)}\right)^\frac{1}{p}.
    \end{align*}
    For $m\in\mathbb{Z}^+$, the Sobolev space $W^{m,p}(\Omega)$ is defined as
    \[
    W^{m,p}(\Omega):=\Big\{f\in L^p(\Omega): D^\alpha f\in L^p(\Omega)\text{ for every }\alpha \in \mathbb{N}^n \text{ with }|\alpha|\leq m\Big\},
    \]
    endowed with the graph norm
    \begin{align*}
        \|f\|_{W^{m,p}(\Omega)}:=\left(\sum_{|\alpha|\leq m}\|D^\alpha f\|_{L^p(\Omega)}^p\right)^{1/p},
    \end{align*}
    where $\alpha=(\alpha_1,\dots,\alpha_n)\in\mathbb{N}^n$, $|\alpha|=\sum_{i=1}^n \alpha_i$, and $D^\alpha f := \frac{\partial^{|\alpha|} f}{\partial x_1^{\alpha_1} \dots \partial x_n^{\alpha_n}}$.
    \item If $s\in(0,1)$, we define the Sobolev-Slobodeckij space $W^{s,p}(\Omega)$ as the space
    \[
    W^{s,p}(\Omega):=\left\{f\in L^p(\Omega):\frac{|f(x)-f(y)|}{|x-y|^{\frac{n}{p}+s}}\in L^p(\Omega\times\Omega)\right\}
    \]
    and we endow $W^{s,p}(\Omega)$ with the graph norm
    \begin{align*}
        \|f\|_{W^{s,p}(\Omega)}&=\left(\|f\|_{L^p(\Omega)}^p+\int_\Omega\int_\Omega\frac{|f(x)-f(y)|^p}{|x-y|^{n+sp}}\mathrm{d}x\mathrm{d}y\right)^\frac{1}{p}\\
        &=:\Big(\|f\|_{L^p(\Omega)}^p+[f]_{W^{s,p}(\Omega)}^p\Big)^\frac{1}{p}.
    \end{align*}
    \item If $s>1$, we denote $s=m+\lambda$, with $m=\lfloor s \rfloor\in\mathbb{Z}$ and $\lambda\in(0,1)$, and set
    \[
    W^{s,p}(\Omega):=\Big\{f\in W^{m,p}(\Omega):D^\alpha f\in W^{\lambda,p}(\Omega)\text{ for every }\alpha\in\mathbb{N}^n: |\alpha|=m\Big\}.
    \]
    As before, we endow $W^{s,p}(\Omega)$ with the norm\footnote{Alternatively, one may define the norm on $W^{s,p}(\Omega)$ using only the Gagliardo seminorm $[D^\alpha f]_{W^{\lambda,p}(\Omega)}$ instead of the full $W^{\lambda,p}$ norm, in order to avoid redundancy, since the $L^p$-norms of derivatives of order $m$ are already included in $\|f\|_{W^{m,p}(\Omega)}$.}
    \begin{align*}
        \|f\|_{W^{s,p}(\Omega)}=\left(\|f\|_{W^{m,p}(\Omega)}^p+\sum_{|\alpha|=m}\|D^\alpha f\|_{W^{\lambda,p}(\Omega)}^p\right)^\frac{1}{p}.
    \end{align*}
    \item If $s=m\in\mathbb{Z}^+$, then, up to equivalent norms,
    \[
    W^{s,p}(\Omega)= W^{m,p}(\Omega).
    \]
\end{itemize}
\end{definition}

Integer-order Sobolev spaces $W^{m,p}(\Omega)$ are not stable neither by real nor complex interpolation.
In the next result we collect some results showing that fractional Sobolev spaces $W^{s,p}(\Omega)$ can always be obtained through real interpolation of integer-order Sobolev families.
The situation is instead more delicate for the complex interpolation of integer-order Sobolev families, whenever $p\neq 2$.\footnote{Given the purpose of this note, the authors have chosen not to introduces Besov spaces and Bessel spaces, which are the correct framework for the link between interpolation and Sobolev spaces.
In the particular cases stated in the Theorem~\ref{thm:fractional_sobolev_interpolation}, they coincide with the Sobolev spaces under consideration.}
We refer the interested reader e.g. to \cite{BL,L} for a proof of the next result, in its most general form.\footnote{We remind that $L^p(\Omega)=W^{0,p}(\Omega)$.}

\begin{theorem}\label{thm:fractional_sobolev_interpolation}
Let $\Omega$ be an open and bounded subset of $\mathbb{R}^n$ of class $\mathbf{C}^1$ and let $p\in[1,+\infty)$.
The following hold:
\begin{itemize}
    \item For every $s\in(0,1)$, up to equivalent norms,
    \[
    W^{s,p}(\Omega)=\Big(L^p(\Omega),W^{1,p}(\Omega)\Big)_{s,p},\quad W^{s,2}(\Omega)=\Big(L^2(\Omega),W^{1,2}(\Omega)\Big)_{[s]}.
    \]
    If $p\neq 2$, then $W^{s,p}(\Omega)\neq\Big(L^p(\Omega),W^{1,p}(\Omega)\Big)_{[s]}$.
    \item More generally, if $s > 1$ and $s \notin \mathbb{Z}^+$, write $s = m + \lambda$ with $m = \lfloor s \rfloor \in \mathbb{Z}^+$ and $\lambda \in (0,1)$. Then, up to equivalent norms,
    \[
    W^{s,p}(\Omega) = \big(W^{m,p}(\Omega), W^{m+1,p}(\Omega)\big)_{\lambda,p},\quad
    W^{s,2}(\Omega) = \big(W^{m,2}(\Omega), W^{m+1,2}(\Omega)\big)_{[\lambda]}.
    \]
    If $p\neq 2$, then $W^{s,p}(\Omega) \neq \big(W^{m,p}(\Omega), W^{m+1,p}(\Omega)\big)_{[\lambda]}$.
\end{itemize}
\end{theorem}
\begin{remark}
For $p\neq 2$, complex interpolation of Sobolev spaces gives rise to Bessel potential spaces $H^{s,p}(\Omega)$ rather than Sobolev-Slobodeckij fractional spaces $W^{s,p}(\Omega)$. 
Similarly, real interpolation at integer indices yields Besov spaces $B^{m}_{p,p}(\Omega)$, which strictly contain $W^{m,p}(\Omega)$ whenever $p\neq 2$. 
These distinctions disappear in the Hilbert case $p=2$, where $H^{s,2}(\Omega)=W^{s,2}(\Omega)$ and $B^{m}_{2,2}(\Omega)=W^{m,2}(\Omega)$.
\end{remark}

As a consequence of Theorem~\ref{thm:inclusions_interpolation}, Theorem~\ref{thm:inclusions_interpolation_complex} and Theorem~\ref{thm:fractional_sobolev_interpolation}, the following embedding results for integer-order and fractional Sobolev spaces follow directly.

\begin{theorem}
Let $\Omega$ be an open subset of $\mathbb{R}^n$ and let $p\in[1,+\infty)$.
The following continuous embeddings hold.
\begin{itemize}
    \item If $s_1,s_2\in(0,1)$ with $s_1\leq s_2$, then
    \[
    W^{s_2,p}(\Omega)\subseteq W^{s_1,p}(\Omega).
    \]
    \item If $\Omega$ is a Lipschitz domain, then for every $s\in(0,1)$ 
    \[
    W^{1,p}(\Omega)\subseteq W^{s,p}(\Omega).
    \]
    \item If $\Omega$ is a Lipschitz domain and $s_1,s_2>1$ satisfy $s_1\leq s_2$, then 
    \[
    W^{s_2,p}(\Omega)\subseteq W^{s_1,p}(\Omega).
    \]
\end{itemize}
\end{theorem}


\section{Mixed local-nonlocal operators}
In the final part of this survey, we consider the so-called \emph{mixed local-nonlocal operators}.
These operators arise from the superposition of a local operator (e.g., the Laplacian $-\Delta u$) and a nonlocal operator (e.g., the fractional Laplacian $(-\Delta)^s u$).
In recent years, there has been growing interest in studying such operators, both for the mathematical challenges they pose and for their relevance in various applications.

The literature on the subject is quite extensive; we refer the interested reader, for instance, to \cite{BDVV,DPLV}.
In \cite{CCMV,CCMV23,MMV}, the authors studied the dependence of mixed local-nonlocal operators on a real parameter (not necessarily positive) $\alpha$.
Let
\[
\mathcal{L}_\alpha u := -\Delta u +\alpha (-\Delta)^s u,
\]
for $\alpha \in \mathbb{R}$ with no a priori restrictions.
We recall that, for fixed $s\in (0,1)$, the fractional Laplacian $(-\Delta)^s u$ can be defined via the Cauchy principal value as
\[
(-\Delta)^s u(x) := C(n,s) \,\mbox{P.V.}\int_{\mathbb{R}^n}\frac{u(x)-u(y)}{|x-y|^{n+2s}}\, dy.
\]
Depending on the value of $\alpha$, we can distinguish the following cases:
\begin{itemize}
    \item When $\alpha=0$, we recover the local Laplacian $-\Delta u$.
    \item When $\alpha >0$, the operator $\mathcal{L}_\alpha$ is still positive-definite. It can be interpreted as the infinitesimal generator of a stochastic process involving both Brownian motion and a pure jump L\'{e}vy process. 
    \item When $-\frac{1}{C}<\alpha<0$, where $C>0$ is the constant of the continuous embedding $W^{1,2}_{0}(\Omega)\subseteq W^{s,2}(\Omega)$, i.e.
    \begin{equation*}
        [u]_{W^{{s,2}}(\Omega)}^2 := \iint_{\mathbb{R}^{2n}}\dfrac{|u(x)-u(y)|^2}{|x-y|^{n+2s}}\, dx dy \leq C \, \int_{\Omega}|\nabla u|^2 \, dx\,,
    \end{equation*}
    the operator remains coercive in this range as well.
    \item Finally, when $\alpha\leq-\tfrac{1}{C}$, $\mathcal{L}_\alpha$ is not positive-definite, the bilinear form naturally associated to it does not induce a scalar product, nor a norm, the variational spectrum may exhibit negative eigenvalues and the maximum principles may fail.
\end{itemize}
Let us show how interpolation theory is useful in such a delicate context. 
Let
\begin{equation*}
\mathbb{H}(\Omega) := \big\{u\in W^{1,2}(\mathbb{R}^n):\,
\text{$u\equiv 0$ a.e.\,on $\mathbb{R}^n\setminus\Omega$}\big\}.
\end{equation*}
If $\Omega$ is, for instance, of class $\mathbf{C}^1$, we identify the space $\mathbb{H}(\Omega)$ with the space $W^{1,2}_0(\Omega)$ as
\begin{equation} \label{eq.identifXWzero}
    u \in W^{1,2}_0(\Omega)\,\,\Longleftrightarrow\,\,
    u\cdot\mathbf{1}_{\Omega}\in\mathbb{H}(\Omega)\,,
\end{equation}
where $\mathbf{1}_\Omega$ is the indicator function of $\Omega$.\footnote{The hypothesis of $\Omega$ of class $\mathbf{C}^1$ can be removed, by working with the closure of regular function.}

From now on, we shall always identify $u\in W^{1,2}_0(\Omega)$ with $\hat{u} := u\cdot\mathbf{1}_\Omega\in\mathbb{H}(\Omega)$.
By the Poincar\'e inequality, the quantity
\[
\|u\|_{\mathbb{H}} :=\left( \int_{\Omega}|\nabla u|^2\, dx\right)^{1/2},\quad u\in\mathbb{H}(\Omega)\,,\]
endows $\mathbb{H}(\Omega)$ with a structure of (real) Hilbert space, which is isometric to $W^{1,2}_0(\Omega)$.
The space $\mathbb{H}(\Omega)$ is separable and reflexive, $\mathbf{C}_0^\infty(\Omega)$ is dense in $\mathbb{H}(\Omega)$ and $\mathbb{H}(\Omega)$ compactly embeds into $L^{2}(\Omega)$ and into 
\[
W^{s,2}_0(\Omega):=\left\{W^{s,2}(\mathbb{R}^n):\,\text{$u\equiv 0$ a.e.\,on $\mathbb{R}^n\setminus\Omega$}\right\}.
\]
Let us consider the following Dirichlet boundary value problem
\begin{equation}\label{eq.EigenvalueProblem}
\begin{cases}
  \mathcal{L}_\alpha u = \lambda u & \textrm{in } \Omega,\\
  u= 0 & \textrm{in } \mathbb{R}^n \setminus \Omega.
\end{cases}
\end{equation}

\begin{definition}
We say that $\lambda \in \mathbb{R}$ is a Dirichlet eigenvalue of $\mathcal{L}_\alpha$ if there exists a weak solution $u\in \mathbb{H}(\Omega)$ of \eqref{eq.EigenvalueProblem} or, equivalently, if
\begin{equation*}
\mathcal{B}_{\alpha}(u,\varphi):=\int_{\Omega}\langle \nabla u, \nabla \varphi\rangle\, dx +\alpha \iint_{\mathbb{R}^{2n}}\dfrac{(u(x)-u(y))(\varphi(x)-\varphi(y))}{|x-y|^{n+2s}}\, dxdy =\lambda \int_{\Omega}u\varphi \, dx
\end{equation*}
for every $\varphi \in \mathbb{H}(\Omega)$.
If such function $u$ exists, we call it eigenfunction corresponding to the eigenvalue $\lambda$.
\end{definition}

The final result of this exposition is a complete description of the (discrete) spectrum of $\mathcal{L}_\alpha$.

\begin{theorem}\label{Thm:spectral}
For every fixed value of $\alpha\in\mathbb{R}$ the following statements hold true:
\begin{itemize}
    \item[($a$)]\, The operator $\mathcal{L}_\alpha$ admits a divergent sequence of eigenvalues $\{\lambda_k\}_{k\in \mathbb{N}}$, bounded from below.
    That is, there exists $C>0$ such that
    \[
    -C < \lambda_1 \leq \lambda_2  \leq \ldots \leq \lambda_k \to +\infty\,, \quad  \textrm{as } k \to +\infty.
    \]
    Moreover, for every $k\in \mathbb{N}$, the eigenvalue $\lambda_{k}$ admits the variational characterization
    \begin{align}\label{eq:lambdaChar}
    \begin{split}
        \lambda_{k} &= \min_{u \in \mathbb{P}_{k},\|u\|_{L^{2}(\Omega)}=1}\left\{\int_{\Omega}|\nabla u|^2 \, dx +\alpha \iint_{\mathbb{R}^{2n}}\dfrac{|u(x)-u(y)|^2}{|x-y|^{n+2s}}\, dxdy \right\}\\
    &=\min_{u \in \mathbb{P}_{k},\|u\|_{L^{2}(\Omega)}=1}\mathcal{B}_{\alpha}(u,u).
    \end{split}
    \end{align}
    Here $\mathbb{P}_{1} := \mathbb{H}(\Omega)$, and, for every $k \geq 2$,
    \begin{equation*}
    \mathbb{P}_{k} := \left\{ u \in \mathbb{H}(\Omega): \mathcal{B}_{\alpha}(u, u_j)=0 \, \textrm{ for every } j = 1, \ldots,k-1\right\},
    \end{equation*}
    with $u_j$ eigenfunction corresponding to the eigenvalue $\lambda_j$.
    \item[($b$)]\, For every $k\in \mathbb{N}$, there exists an eigenfunction $u_{k} \in \mathbb{H}(\Omega)$ corresponding to $\lambda_{k}$, which realizes the minimum in \eqref{eq:lambdaChar}.
    \item[($c$)]\, The sequence $\{u_k\}_{k\in \mathbb{N}}$ of eigenfunctions is an orthonormal basis of $L^{2}(\Omega)$, and orthogonal also with respect to the bilinear form $\mathcal{B}_\alpha$.
    \item[($d$)]\, Each eigenvalue $\lambda_k$ has finite multiplicity.
\end{itemize}
\end{theorem}

The proof of Theorem~\ref{Thm:spectral} can be found in \cite{CCMV,MMV} and a sketch of it is reported below.
The key point of the proof is the following interpolation inequality, which follows from \eqref{ext:Brezis_Miro}, and whose proof can be found in \cite{interpolation}.

\begin{theorem}\label{thm:Brezis}
Let $\Omega$ be a Lipschitz bounded domain in $\mathbb{R}^n$.
There exists a positive constant $C$ such that
\begin{equation*}
    [u]_{W^{{s,2}}(\Omega)}^2\leq C\|u\|^{2(1-s)}_{L^2(\Omega)}\|u\|^{2s}_{W^{1,2}(\Omega)}
\end{equation*}
for every $u\in\mathbb{H}(\Omega)$. 
\end{theorem}
The optimal constant $C$ is explicitly computed with the help of Fourier transform; see e.g., \cite{CCMV,CCMV23}.
    
\begin{proof}[Sketch of the proof of Theorem~\ref{Thm:spectral}]
We limit this brief exposition to demonstrating the existence of eigenvalues for the operator $\mathcal{L}_\alpha$.
We distinguish the following scenarios, according with the value of $\alpha$.
\begin{itemize}
    \item If $\alpha >-\tfrac{1}{C}$, $C>0$ being the constant of the embedding $W^{1,2}_0(\Omega)\subseteq W^{s,2}(\Omega)$, the conclusions follows by the classical spectral theorem.
    \item In what follows we assume $\alpha \leq -\tfrac{1}{C}$.
\end{itemize}
By Theorem~\ref{thm:Brezis} and the Young inequality, with exponents $\frac{1}{s}$ and $\frac{1}{1-s}$, for every $\varepsilon\in\mathbb{R}^+$ there exist positive constants $c_{\varepsilon},c_1,c_2$, depending only on $s$ and $\varepsilon$, such that
\begin{equation*}
\begin{split}
    |\alpha|[u]^2_s&\leq |\alpha|C\left((1-s)c_\varepsilon\|u\|^2_{L^2(\Omega)}+s\varepsilon\|u\|^2_{W^{1,2}(\Omega)}\right)\\
    &=|\alpha|c_1\varepsilon\|u\|^2_{\mathbb{H}(\Omega)}+|\alpha|c_2\|u\|^2_{L^2(\Omega)}
\end{split}
\end{equation*}
for every $u\in\mathbb{H}(\Omega)$.
By choosing $\varepsilon=\frac{1}{2c_1|\alpha|}$, we get that
\begin{equation*}
    |\alpha|[u]^2_s\leq\frac{1}{2}\|u\|^2_{\mathbb{H}(\Omega)}+\gamma\|u\|^2_{L^2(\Omega)}
\end{equation*}
for every $u\in\mathbb{H}(\Omega)$, where $\gamma$ only depends on $s,\varepsilon$ and $\alpha$.
Therefore,
\begin{equation*}
\mathcal{B}_\alpha(u,u)+\gamma\|u\|^2_{L^2(\Omega)}\geq\frac{1}{2}\|u\|^2_{\mathbb{H}(\Omega)}\quad\text{for every }u\in\mathbb{H}(\Omega).
    \end{equation*}
Then, the bilinear form $\mathcal{B}_\alpha^\gamma \colon \mathbb{H}(\Omega)\times \mathbb{H}(\Omega)\to \mathbb{R}$, defined as
\begin{equation*}
\mathcal{B}_\alpha^\gamma (u,v):=\mathcal{B}_\alpha(u,v)+\gamma (u,v)_{L^2(\Omega)}\quad\text{for all $u,v\in \mathbb{H}(\Omega)$},
\end{equation*}
is symmetric, continuous, and coercive.

The associated linear and continuous operator $\mathcal{L}_\alpha^\gamma$ is related with $\mathcal{L}_\alpha$ by 
\[
\mathcal{L}_\alpha^\gamma u=\mathcal{L}_\alpha u+\gamma u\quad\text{ for every }u\in \mathbb{H}(\Omega).
\]
By Lax-Milgram Lemma, $\mathcal{L}_\alpha^\gamma$ is invertible and the inverse operator $(\mathcal{L}_\alpha^\gamma)^{-1}$ is linear and continuous.
We also note that $\lambda\in\mathbb{R}$ is an eigenvalue of $\mathcal{L}_\alpha$ in $\mathbb{H}(\Omega)$, with eigenvector $u\in \mathbb{H}(\Omega)\setminus\{0\}$, if and only $\frac{1}{\lambda+\gamma}$ is an eigenvalue of $R^\gamma:=(\mathcal{L}_\alpha^\gamma)^{-1}$ in $L^2(\Omega)$, with eigenvector $u\in L^2(\Omega)\setminus\{0\}$.
The operator $R^\gamma\colon L^2(\Omega)\to L^2(\Omega)$ is linear, continuous, and compact, since the embedding $\mathbb{H}(\Omega)\subset L^2(\Omega)$ is compact. Moreover, $R^\gamma$ is injective and self-adjoint, being $\mathcal{B}_\alpha^\gamma $ symmetric, and
\begin{equation*}
(R^\gamma h,h)_{L^2(\Omega)}=\mathcal{B}_\alpha^\gamma (R^\gamma h,R^\gamma h)\ge 0\quad\text{for all $h\in L^2(\Omega)$}.
\end{equation*}
Therefore, there exists a decreasing sequence $(\mu_k)_k$ of eigenvalues of $R^\gamma$ with $\mu_k>0$ and $\mu_k\to 0$ as $k\to\infty$. Moreover, every $\mu_k$ has finite multiplicity and there exists a orthonormal basis $(e_k)_k$ of $L^2(\Omega)$ given by eigenvectors of $R^\gamma$ associated to $\mu_k$. Hence, if we consider
\begin{equation*}
\lambda_k:=\frac{1}{\mu_k}-\gamma\quad\text{for all $k\in\mathbb N$},
\end{equation*}
we get that $(\lambda_k)_k$ is an increasing sequence of eigenvalues of $\mathcal{L}_\alpha$ in $\mathbb{H}(\Omega)$.

The conclusions follow again from the standard proof of the spectral theorem.
\end{proof}

\noindent{\it \textbf{Acknowledgements.}} The author wishes to thank Augusto Visintin, who introduced him to the topic, and Maicol Caponi, Nicolò Cangiotti, Enrico Valdinoci and Enzo Vitillaro for many useful suggestions.

\medskip

\noindent{\it \textbf{Funding information.}} The author is a member of GNAMPA of the Istituto Nazionale di Alta Matematica (INdAM), and acknowledges the support of the INdAM - GNAMPA 2025 Project ``Metodi variazionali per problemi dipendenti da operatori frazionari isotropi e anisotropi'' (Grant Code: CUP\_E5324001950001).
He is also supported by the Spanish State Research Agency, through the Severo Ochoa and María de Maeztu Program for Centers and Units of Excellence in R\&D (CEX2020-001084-M), by MCIN/AEI/10.13039/501100011033 (PID2021-123903NB-I00), and by Generalitat de Catalunya (2021-SGR-00087).

\end{document}